\documentclass[12pt]{article}
\author{Thomas Morrill\footnote{Supported by Australian Research Council Discovery Project DP160100932 }\\ School of Physical, Environmental and Mathematical Sciences\\ The University of New South Wales Canberra, Australia\\t.morrill@adfa.edu.au\and  Tim Trudgian$^{*}$\footnote{Supported by Australian Research Council Future Fellowship FT160100094.} \\
School of Physical, Environmental and Mathematical Sciences\\ The University of New South Wales Canberra, Australia \\
  t.trudgian@adfa.edu.au}
  \title{An elementary bound on Siegel zeroes}
\usepackage{indentfirst}
\usepackage{url}
\usepackage{enumerate}
\usepackage{amsthm}
\usepackage{amsmath}
\usepackage{comment}
\usepackage{fullpage}
\usepackage{amssymb}
\usepackage{booktabs}
\usepackage{siunitx}
\usepackage{xcolor}
\usepackage{graphicx}

\newtheorem{thm}{Theorem}

\newtheorem{Lem}{Lemma}

\begin{document}
\maketitle
\begin{abstract}
\noindent 
We consider Dirichlet $L$-functions $L(s, \chi)$ where $\chi$ is a real, non-principal character modulo $q$. Using Pintz's refinement of Page's theorem, we prove that for $q\geq 3$ the function $L(s, \chi)$ has at most one real zero $\beta$ with $1- 1.011/\log q < \beta < 1$.
\end{abstract}
\section{Introduction}
\noindent
Let $\chi(n)$ be a real, non-principal Dirichlet character to the modulus $q$ and let $L(s, \chi)$ be the associated Dirichlet $L$-function, where $s=\sigma + it$. It is known \cite[pp.\ 93--95]{DP} that $L(s, \chi)$ has at most one zero with real part larger than $1 - (A\log \max\{q, q|t|\})^{-1}$. Such an exceptional, or Siegel, zero must lie on the real axis. 
\noindent
A classical result of Page \cite{Page} is that, given a single character $\chi$ modulo $q$, there can be at most one exceptional zero `close' to unity.
\begin{thm}[Page]\label{jack}
If $\chi$ mod $q$ is a real, non-principal character, and if $\beta_{1}$ and $\beta_{2}$ are real zeroes of $L(s, \chi)$, then there is a positive constant $c_1$ such that
$$\min(\beta_{1}, \beta_{2}) \leq 1 - \frac{c_1}{ \log q}.$$
\end{thm}

McCurley \cite{McCurley} and Kadiri \cite{Kadiri} have given values for $c_1$: the best is $c = 0.909$ by Kadiri. Kadiri's method, similar to that in papers on the zero-free region of $L$-functions (see, e.g.\ \cite{Kadiri2} and \cite{MT}), uses a special nonnegative trigonometric polynomial, the calculus of variations, and an analysis of the distribution of the imaginary parts of zeroes of $L(s, \chi)$.

Pintz \cite[Thm 2]{Pintz} revisits Page's method, which is more elementary. Using the P\'{o}lya--Vinogradov inequality, Pintz is able to prove that $c_1 = 2 + o(1)$ as $q\rightarrow\infty$. Indeed, he notes that one can use the Burgess bounds on character sums to improve this to $c_1 = 4 + o(1)$. In making these results explicit there will be some loss in the size of $c_1$. We aim to minimise this loss by using the best off-the-shelf explicit estimates. Our main result is the following.

\begin{thm}\label{jill}
If $\chi$ mod $q$ is a real, non-principal character, with $q\geq 3$, and if $\beta_{1}$ and $\beta_{2}$ are real zeroes of $L(s, \chi)$, then
$$\min(\beta_{1}, \beta_{2}) \leq 1 - \frac{1.011}{ \log q}.$$
\end{thm}

Throughout the course of the paper we take take $\chi$ to be a primitive character. This is no great obstacle, since, as noted by Pintz \cite[p.\ 164]{Pintz}, if $\chi$ modulo $q$ is induced by a primitive character $\chi'$ modulo $q'$, then if $L(\beta, \chi)= 0$ we have $L(\beta, \chi') = 0$. Since $q \leq q'$, we can therefore extend the result to that in Theorem \ref{jill}.

The rest of this paper is organised as follows. In \S\ref{snooker} we collect some preliminary results. We use these, in \S\ref{nine} to refine Pintz's result for finite ranges of $q$. We detail, in \S\ref{mine}, our computations. These prove Theorem \ref{jill} for finite ranges of $q$; we then use Pintz's original argument for large $q$. Finally, in \S\ref{otello} we outline some potential improvements to our results.

Throughout the paper $\vartheta$ will denote a complex number with modulus at most unity.

\section{Preliminary lemmas}\label{snooker}
In this section we collect some results from the literature. We first note that we need not concern ourselves with small values of $q$. Watkins \cite{Watkins} showed that there are no Siegel zeroes for $L(s, \chi)$ where $\chi$ is odd, and $q\leq 3\cdot 10^{8}$; Platt \cite{Platt} reached the same conclusion for $\chi$ even and $q\leq 4\cdot 10^{5}$. 

We wish to record an explicit version of the P\'{o}lya--Vinogradov inequality.
\begin{Lem}[Frolenkov and Soundararajan \cite{Sound}] \label{egg}
Let $\chi$ be a primitive character with parity $\chi(-1) = (-1)^i$.
We have for all $q\geq 1200$, that
\begin{equation}\label{chalk}
	\left| \sum_{n=M+1}^{M+N} \chi(n) \right| \leq E_i \sqrt{q} \log q + \sqrt{q},
\end{equation}
where $E_0 = {2}/{\pi^2}$, $E_1 = {1}/{2 \pi}$.
\end{Lem}
We note that, by a remark of Pomerance \cite{Pomerance}, we can divide the right-side of (\ref{chalk}) by two if $M=0$ and $\chi$ is even.
Define 
\begin{align}\label{macbet}
	A (q_0)
	= A_i(q_0)
	:= \begin{cases}
		{1}/{\pi^2} + 1/(2\log q_{0}), &  \text{if} \; i = 0\\
		{1}/{2 \pi} + 1/\log q_{0}, &  \text{if} \; i = 1,
	\end{cases}
\end{align}
so that, for $q\geq q_{0} \geq 1200$,
\begin{align*} \label{dust}
	\left| \sum_{n=M+1}^{M+N} \chi(n) \right|
	\leq A(q_0) \sqrt{q} \log q.
\end{align*}

We shall make use of the Euler--Maclaurin summation formula --- see \cite[Thm 2.19]{Murty}.
\begin{Lem}[Euler--Maclaurin summation]\label{EML}
Let $k$ be a nonnegative integer and $f(x)$ be $(k+1)$ times differentiable on the interval $[a, b]$. Then
\begin{equation*}\label{Ems}
\begin{split}
\sum_{a<n\leq b} f(n) &= \int_{a}^{b} f(t)\, dt + \sum_{r=0}^{k} \frac{(-1)^{r+1}}{(r+1)!} \left(f^{(r)}(b) - f^{(r)}(a)\right) B_{r+1}\\
&+ \frac{(-1)^{k}}{(k+1)!} \int_{a}^{b} B_{k+1}(x) f^{(k+1)}(x)\, dx,
\end{split}
\end{equation*}
where $B_{j}(x)$ is the $j$th periodic Bernoulli polynomial and $B_{j} = B_{j}(0)$.
\end{Lem}
Pintz takes $k=0$ and examines sums of $n^{-\alpha}$ and $(\log n)n^{-\alpha}$. We shall require more precision for our calculations. Choosing $k=2$ we have, for $0<\alpha<1$,
\begin{equation}\label{quince}
\sum_{n\leq x} n^{-\alpha} = C_{1}(\alpha) + \frac{1}{1-\alpha} \left(x^{1-\alpha}-1\right) + \frac{1}{2}x^{-\alpha} - \frac{\alpha}{12} x^{-\alpha -1} + \vartheta \frac{\alpha(\alpha +1)}{72 \sqrt{3}}x^{-\alpha -2}, 
\end{equation}
where 
\begin{equation*}\label{paste}
C_{1}(\alpha) = \frac{1}{2} + \frac{\alpha}{12}  + \frac{\alpha(\alpha+1)(\alpha+2)}{6} \int_{1}^{\infty} \frac{\{t\}^{3} - \frac{3}{2}\{t\}^{2} + \frac{1}{2} \{t\}}{t^{\alpha + 3}}\, dt.
\end{equation*}
Similarly, with $k=2$ we have for $0<\alpha <1$
\begin{equation*}\label{boa}
\begin{split}
\sum_{n\leq x} (\log n)n^{-\alpha} =& C_{2}(\alpha) + \frac{x^{1-\alpha}\log x}{1-\alpha} + \frac{1}{(1-\alpha)^{2}}\left(1 - x^{1-\alpha}\right) + \frac{1}{2} x^{-\alpha} \log x \\
&+ \vartheta x^{-\alpha -1}\left\{   \alpha \log x\left( 1 + \frac{1}{72 \sqrt{3}}\right) + 1 + \frac{3\alpha + 1}{72\sqrt{3}(\alpha +1)} \right\},
\end{split}
\end{equation*}
where
\begin{equation*}\label{constrictor}
C_{2}(\alpha) = 1  - \frac{1}{12} + \frac{1}{6} \int_{1}^{\infty}\frac{(\alpha(\alpha + 1)\log t - 1-2\alpha)(\{t\}^{3} - \frac{3}{2}\{t\}^{2} + \frac{1}{2} \{t\})}{t^{\alpha +2}}\, dt.
\end{equation*}

The class-number formula allows one to show \cite[p.\ 95]{DP} that $L(1, \chi)\gg q^{-1/2}$. We require an explicit version of this as given by Bennett, Martin, O'Bryant and Rechnitzer \cite[Lem.\ 6.3 and A.9]{Bennett}.
\begin{Lem}[Bennett et al.]\label{trovatore} 
If $\chi$ modulo $q$ is a real primitive character, then
\begin{equation*}\label{harp}
	L(1, \chi) \geq
	\begin{cases}
		79.2q^{-1/2}, & \mbox{if } 4\cdot 10^{5} \leq q \leq 10^{7} \\ 
		12 q^{-1/2}, & \mbox{if }  q> 10^{7}.
	\end{cases}
\end{equation*}
\end{Lem}

Consider
\begin{equation}\label{busy}
\sum_{n\leq x} \frac{g(n)\log n}{n^{1-\tau}}, \quad g(n) = \sum_{d|n} \chi(d),
\end{equation}
where $0<\tau<1$. As in Pintz, we note that
$$g(n) = \prod_{p^{e}||n}\left(1 + \chi(p) + \chi(p^{2}) + \cdots + \chi(p^{e})\right).$$
Since $\chi(n)$ is absolutely multiplicative we have that $g(m^{2}) \geq 1$.
Pintz uses this to show that the first sum in (\ref{busy}) exceeds $(\log 4)/4$ for $x\geq 4$.
We improve this, using partial summation.
\begin{Lem}\label{eagle}
\begin{equation}\label{james}
\sum_{n\leq x} \frac{g(n)\log n}{n^{1-\tau}} \geq -2 \zeta'(2-2\tau)  - f(2\tau -1, x), 
\end{equation}
where
\begin{equation*}\label{suit}
f(\alpha, x) = \frac{2(x^{1/2} -1)^{\alpha} (1 - \alpha \log(x^{1/2} -1))}{\alpha^{2}}, \quad \left(\log(x^{1/2} -1) > \frac{1}{2\tau -1}\right).
\end{equation*}
\end{Lem}
We note that we shall only apply (\ref{james}) for finite values of $x$, and, as such, we can avoid the usual irritation about bounding terms such as $x^{1/2} -1$ from below.

Finally, by partial summation we have, for any continuously differentiable function $h(z)$,
\begin{align*} \label{Saturn}
	\sum_{x<n\leq y} \chi(n) h(n) =
	h(y) \sum_{1\leq n \leq y} \chi(n)
	- h(x) \sum_{1 \leq n\leq x} \chi(n)
	- \int_{x}^{y} h'(t) \sum_{1\leq n\leq t} \chi(n)\, dt.
\end{align*}
Suppose that $h$ is a positive decreasing function with $\lim_{z \to \infty} h(z) = 0$.
For our purposes, $h \in \{(\log z) / z^{1 - \tau}, \; 1 / z^{1 - \tau}, \; 1/z \}$.
For any $x$ we have $\sum_{n\leq x} \chi(n) \ll \sqrt{q} \log q$, by the P\'{o}lya--Vinogradov inequality, whence
\begin{equation}\label{baulk}
	\left|\sum_{n > x} \chi(n) h(n)\right|
	\leq h(x) \left|\sum_{1 \leq n \leq x} \chi(n)\right|
	+ \int_{x}^{\infty} |h'(t)|\left| \sum_{1 \leq n \leq t} \chi(n)\right|\, dt.
\end{equation}
Depending on the size of $z$ we shall use a mixture of the trivial bound and the P\'{o}lya--Vinodgradov inequality in (\ref{baulk}). Doing this, and integrating by parts, yields

\begin{equation}\label{hawk}
	\left|\sum_{n>z} \chi(n) h(n)\right|
	\leq \begin{cases}
		2 z h(z) + \int_{z}^{A (q_0) \sqrt{q} \log q} h(t)\, dt, & \textrm{if}\; z\leq A(q_0)\sqrt{q}\log q\\
		2 A (q_0) \sqrt{q} (\log q) h(z), & \textrm{if}\; z \geq A(q_0)\sqrt{q}\log q,
	\end{cases}
\end{equation}
where $A(q_{0})$ is defined in (\ref{macbet}).

\section{Outline of proof}\label{nine}
Following Pintz, we apply our Lemma \ref{EML} to obtain
\begin{equation}\label{cello}
\sum_{n\leq x} \frac{g(n)\log n}{n^{1-\tau}} = T_{1} + T_{2} - T_{3} - \frac{(1 - \tau)}{12} x^{\tau-2} \sum_{d\leq x} \chi(d) d \log d + \frac{1}{2} x^{\tau -1} \log x \sum_{d\leq x}\chi(d) + \vartheta W,
\end{equation}
where
\begin{align*}
	T_1
	&= K_1 \sum_{d \leq x} \chi(d) \frac{\log d}{d^{1 - \tau}} \\
		T_2
		&= K_2 \sum_{d \leq x} \chi(d) \frac{1}{d^{1 - \tau}} \\
			T_3
			&= \frac{x}{\tau} \bigg(\frac{1}{\tau} - \log x\bigg) \sum_{d \leq x} \chi(d) \frac{1}{d}.
\end{align*}
Here, $K_1$ and $K_2$ are positive constants depending only on $\tau$, and
\begin{equation*}\label{fret}
\begin{split}
	W = \;
	& x^{\tau} \log x \frac{(1-\tau)(2-\tau)}{432 \sqrt{3}} \left( 2 + x_{0}^{-1}\right)\left( 1 + x_{0}^{-1}\right) \\
	& + \frac{x^{\tau}}{2} \left( 1 + x_{0}\right)^{-1}\left((1-\tau)   \left( 1 + \frac{1}{72 \sqrt{3}}\right)\log x + 1 + \frac{4-3\tau}{72 \sqrt{3}(2-\tau)}\right).
\end{split}
\end{equation*}
Note that Pintz introduces the variable $z \leq x$ and splits the left-side of \eqref{cello}  to bound what would otherwise be $x^{1 + \tau} \log x$ in $(3.5)$ of \cite{Pintz}.
This is rendered unnecessary by using \eqref{quince}.

We bound the character sums in (\ref{cello}) by using \eqref{hawk}. This produces
\begin{align*}\label{harrier}
	\sum_{n\leq x} \frac{g(n)\log n}{n^{1-\tau}}
	 = K_1 L'(1-\tau)
	+ K_2 L(1-\tau)
	 - \frac{x}{\tau} \bigg(\frac{1}{\tau} - \log x\bigg)  L(1)
	+ E(q, \tau, x),
\end{align*}
where $E$ is an unwieldy, though easily computed, error term.
Assume, that there are two zeroes of $L(s, \chi)$ with $s$ in the interval $(1-c/\log q, 1)$. Hence there is a value of $\tau\in(0, c/\log q)$ for which $L(1-\tau) \leq 0$ and $L'(1-\tau) = 0$. We therefore have
\begin{equation}\label{charles}
	\sum_{n\leq x} \frac{g(n)\log n}{n^{1-\tau}}
	\leq E(q, \tau, x)
	- \left( \frac{1}{\tau} - \log x\right) \frac{B x^{\tau}}{\tau \sqrt{q}},
	\quad (\tau^{-1} \geq \log x),
\end{equation}
where $B$ is either $79.2$ or $12$, according to Lemma \ref{trovatore}.
We now invoke Lemma \ref{eagle}, which provides us with a contradiction if $F= F(q, \tau, x)<0$, where
\begin{equation*}\label{oliver}
	F
	= E(q, \tau, x)
	- \left( \frac{1}{\tau} - \log x\right) \frac{B x^{\tau}}{\tau \sqrt{q}}
	+ 2 \zeta'(2-2\tau)
	+ 2 \left( \frac{ 1 + (\frac{1}{2} - \tau)\log x_{0}}{x_{0}^{1/2 -\tau}(1-2\tau)^{2}}\right)
	< 0, 
\end{equation*}
subject to $x_{0} \leq x \leq \exp(\tau^{-1})$.

\subsection{Algorithm} \label{mine}

We wish to calculate the best constant $c$ in Theorem \ref{jill} on some range $q_{0} \leq q \leq q_{1}$.
We calculate an upper bound $F^\star$ so that $F \leq F^\star$ for all $q \in [q_0, q_1]$.
For a fixed $c$, it suffices to find $0 < x^\star \leq \exp(\tau^{-1})$ so that $F^\star(q, c/\log q, x) < 0$.
The algorithm calculates $F^\star$ at test points $x^\star$; if an $x^\star$ is found so that $F^\star < 0$, the algorithm increments $c$ and restarts the search.
When no admissible $x^\star$ values are found, the algorithm terminates and returns the last known $c$ and $x^\star$ values for which Theorem \ref{jill} is true.

The algorithm is run separately for even and odd characters, using the P\'{o}lya--Vinogradov bounds in Lemma~\ref{egg}.
The results of this computation are given in Table \ref{gloves}.
This computation was run in Sage on a $2.9$ GHz processor.
The code is available at \url{https://github.com/tsmorrill/Pintz}.

\begin{table}
\begin{center}
\begin{tabular}{ll|ll|ll}
	\hline
	\hline
	$q_0$ & $q_1$ & $c_{\textrm{even}}$ & $x_{\textrm{even}}$ & $c_{\textrm{odd}}$ & $x_{\textrm{odd}}$ \\
	\hline
	${4} \cdot 10^{5}$ & ${7} \cdot 10^{5}$ & \num{1.011} & $10^{5.54}$ & -- & -- \\
	${7} \cdot 10^{5}$ & $ 10^{6}$ & \num{1.017} & $10^{5.73}$ & -- & -- \\
	$ 10^{6}$ & ${1.7} \cdot 10^{6}$ & \num{1.020} & $10^{5.88}$ & -- & -- \\
	${1.7} \cdot 10^{6}$ & ${3.1} \cdot 10^{6}$ & \num{1.025} & $10^{6.08}$ & -- & -- \\
	${3.1} \cdot 10^{6}$ & ${6.1} \cdot 10^{6}$ & \num{1.030} & $10^{6.30}$ & -- & -- \\
	${6.1} \cdot 10^{6}$ & ${1.3} \cdot 10^{7}$ & \num{1.036} & $10^{6.54}$ & -- & -- \\
	${1.3} \cdot 10^{7}$ & ${2.4} \cdot 10^{7}$ & \num{1.041} & $10^{6.83}$ & -- & -- \\
	${2.4} \cdot 10^{7}$ & ${5.4} \cdot 10^{7}$ & \num{1.044} & $10^{7.06}$ & -- & -- \\
	${5.4} \cdot 10^{7}$ & ${1.5} \cdot 10^{8}$ & \num{1.051} & $10^{7.35}$ & -- & -- \\
	${1.5} \cdot 10^{8}$ & ${6.2} \cdot 10^{8}$ & \num{1.055} & $10^{7.75}$ & \num{1.021} & $10^{8.00}$ \\
	${6.2} \cdot 10^{8}$ & ${4.4} \cdot 10^{9}$ & \num{1.060} & $10^{8.29}$ & \num{1.029} & $10^{8.54}$ \\
	${4.4} \cdot 10^{9}$ & ${6.4} \cdot 10^{10}$ & \num{1.070} & $10^{9.01}$ & \num{1.041} & $10^{9.26}$ \\
	${6.4} \cdot 10^{10}$ & ${2.7} \cdot 10^{12}$ & \num{1.080} & $10^{10.00}$ & \num{1.055} & $10^{10.24}$ \\
	${2.7} \cdot 10^{12}$ & ${6.2} \cdot 10^{14}$ & \num{1.090} & $10^{11.40}$ & \num{1.069} & $10^{11.63}$ \\
	${6.2} \cdot 10^{14}$ & ${2.1} \cdot 10^{18}$ & \num{1.101} & $10^{13.43}$ & \num{1.082} & $10^{13.66}$ \\
	${2.1} \cdot 10^{18}$ & ${4.4} \cdot 10^{21}$ & \num{1.200} & $10^{15.26}$ & \num{1.182} & $10^{15.50}$ \\
	${4.4} \cdot 10^{21}$ & ${1.5} \cdot 10^{24}$ & \num{1.300} & $10^{16.64}$ & \num{1.283} & $10^{16.86}$ \\
	${1.5} \cdot 10^{24}$ & ${3.1} \cdot 10^{26}$ & \num{1.350} & $10^{17.90}$ & \num{1.334} & $10^{18.12}$ \\
	${3.1} \cdot 10^{26}$ & ${2.4} \cdot 10^{28}$ & \num{1.400} & $10^{18.92}$ & \num{1.383} & $10^{19.15}$ \\
	${2.4} \cdot 10^{28}$ & ${1.5} \cdot 10^{30}$ & \num{1.425} & $10^{19.91}$ & \num{1.411} & $10^{20.11}$ \\
	${1.5} \cdot 10^{30}$ & $ 10^{32}$ & \num{1.445} & $10^{20.89}$ & \num{1.429} & $10^{21.11}$ \\
	$ 10^{32}$ & ${9.1} \cdot 10^{32}$ & \num{1.495} & $10^{21.40}$ & \num{1.480} & $10^{21.61}$ \\
	\hline
	\hline
\end{tabular}
\end{center}
\caption{Values of $c$ and $x$ so that $F(q, \tau, x) < 0$ for $q \in [q_0, q_1]$. \label{gloves}
}
\end{table}

\subsection{Large moduli} \label{mutiny}

For $q$ outside Table \ref{gloves}, we apply Lemma $3$ of \cite{Pintz} with $x = A\sqrt{q}(\log q)/\tau^{8}$, where $A = A(q_{0})$ as defined in (\ref{macbet}).   
Suppose $L(s, \chi)$ has two zeroes in the interval $(1-c/\log q, 1)$.
Pintz considers a variation of (\ref{charles}): to obtain a contradiction he requires
\begin{equation} \label{calypso}
	\frac{6ec}{\log q}
	< \frac{\log 4}{4}\qquad \textrm{and} \qquad \frac{1}{\tau} - \log x >0.
\end{equation}
For $c = 1.011$, a quick check shows that we need $q\geq 4.6 \cdot 10^{20}$ for the first inequality in (\ref{calypso}) to hold. Consider the second inequality: for $\tau \leq c/\log q$, we have that
\begin{align} \label{fruitcakes}
	1/\tau - \log x
	= 1/\tau - \log A - \log(\sqrt{q}\log q)- 8 \log 1/\tau.
\end{align}

If we have $\tau < 1/8$, then \eqref{fruitcakes} is decreasing -- fortunately, this is implied by $6e c/{\log q} < (\log 4)/{4}$.
Therefore, we have
\begin{align} \label{cheeseburger}
	1/\tau - \log x
	\geq \left(\frac{1}{c} - \frac{1}{2}\right)\log q - \log A - 9 \log\log q + 8 \log c,
\end{align}
for  $q\geq 4.6 \cdot 10^{20}$.
To ensure that the right-side of (\ref{cheeseburger}) is positive,   we need $q > 2.6 \cdot 10^{32}$ for even characters, and $q > 9.1 \cdot 10^{32}$ for odd characters.
Thus, Theorem \ref{jill} holds for all $q > 9.1 \cdot 10^{32}$. This, along with Table \ref{gloves} completes the proof.

We note that the argument leading to the contradiction could be improved by replacing the $(\log 4)/4$ bound with the result from Lemma \ref{eagle}. However, the more difficult inequality to satisfy is \eqref{cheeseburger}, and retaining the $(\log 4)/4$ eases the computation.

\section{Conclusions}\label{otello}
Our result can be improved at several places. A marginally better constant in the P\'{o}lya--Vinogradov inequality in Lemma \ref{egg} gives little overall improvement. Similarly, if one extended the computation done by Bennett et al.\ \cite{Bennett}, and dealt with some small values of $q$ directly, one may improve slightly on the lower bounds on $L(1, \chi)$ in Lemma \ref{trovatore}.

More importantly, the small values of $q$ we are forced to consider impede our calculation of $c$. If Platt's result were extended to show that there are no Siegel zeroes for $L(s, \chi)$ for $\chi$ even and $q \leq Q$ where $Q> 4\cdot 10^{5}$, then Theorem \ref{jill} may be improved according to Table \ref{gloves}. Note that odd characters  above $3\cdot 10^{8}$ must also be dealt with to improve $c \geq 1.02$.

We have essentially `lost half' of Pintz's $c_{1} = 2 + o(1)$ result in obtaining our Theorem \ref{jill}. This gives  hope to using the Burgess bounds (asymptotically giving $c_{1} = 4 + o(1)$) to improve further on our results. Explicit versions of the Burgess bounds are available (for example, see \cite{Trevino}). One could splice these results with the trivial and P\'olya--Vinogradov estimates.



\begin{thebibliography}{1}
\bibitem{Bennett} M.A. Bennett, G. Martin, K. O'Bryant, and A. Rechnitzer.
\newblock Explicit bounds for primes in arithmetic progressions. To appear in {\em Illinois J. Math}. Preprint available at
\newblock {\em arXiv: 1802.00085}, 2018.

\bibitem{DP} H. Davenport. \newblock {\em Multiplicative Number Theory}. Second edition. Graduate Texts in Mathematics, 74. Springer-Verlag, New York-Berlin, 1980.

\bibitem{Sound}
D.A. Frolenkov and K. Soundararajan.
\newblock{A generalization of the P\'{o}lya--Vinogradov inequality}.
\newblock{\em Ramanujan J.} 31(3), 271--279 (2013).



\bibitem{Kadiri}
H. Kadiri.
\newblock{An explicit zero-free region for Dirichlet $L$-functions}.
\newblock{\em Preprint}.

\bibitem{Kadiri2}
H. Kadiri.
\newblock{Explicit zero-free regions for Dirichlet $L$-functions.}
\newblock{\em Mathematika} 64(2), 445--474 (2018).


\bibitem{McCurley}
K.S. McCurley.
\newblock{Explicit zero-free regions for Dirichlet $L$-functions.}
\newblock{\em J. Number Theory} 19, 7--32 (1984).

\bibitem{MT}
M. J. Mossinghoff and T. S. Trudgian.
\newblock Nonnegative trigonometric polynomials and a zero-free region for the Riemann zeta-function.
\newblock {\em J. Number Theory}, 157, 329--349 (2015).

\bibitem{Murty}
M.R. Murty.
\newblock{\em Problems in Analytic Number Theory}. Second edition. Graduate Texts in Mathematics, 206. Reachings in Mathematics. Springer, New York, 2008.

\bibitem{Page}
A. Page.
\newblock{On the number of primes in an arithmetic progression.}
\newblock{\em Proc. London Math. Soc.} 39, 116--141 (1935).

\bibitem{Pintz}
J. Pintz.
\newblock Elementary methods in the theory of $L$-functions, V. The Theorems of Landau and Page.
\newblock {\em Acta Arith.}, 32,  163--171 (1977).

\bibitem{Platt}
D. Platt.
\newblock{Numerical computations concerning the GRH.} 
\newblock{\em Math. Comp.} 85, 3009--3027 (2015).

\bibitem{Pomerance}
C. Pomerance .
\newblock{Remarks on the P\'{o}lya--Vinogradov inequality}.
\newblock{\em Integers} 11A, A19, 11pp.\ (2009).

\bibitem{Trevino}
E.~Trevi\~{n}o.
\newblock The {B}urgess inequality and the least $k$-th power non-residue.
\newblock {\em Int. J. Number Theory}, 11(5), 1653--1678 (2015).

\bibitem{Watkins}
M. Watkins.
\newblock{Real zeros of real odd Dirichlet $L$-functions.}
\newblock{\em Math. Comp.} 73, 415--423 (2004).


\end{thebibliography}
 \end{document}